\documentclass[11pt,twoside]{article}

\usepackage{amsthm,amsmath,amssymb}

\begin{document}
\newcommand{\ol }{\overline}
\newcommand{\ul }{\underline }
\newcommand{\ra }{\rightarrow }
\newcommand{\lra }{\longrightarrow }
\newcommand{\ga }{\gamma }

\title{{Some Functorial Properties of Nilpotent Multipliers}}
\author{Behrooz Mashayekhy and Mahboobeh Alizadeh Sanati \\ Department of
Mathematics,\\ Center of Excellence in Analysis on Algebraic
Structures,\\
Ferdowsi University of Mashhad,\\ P.O.Box 1159-91775, Mashhad, Iran. \\
{\small E-mail: mashaf@math.um.ac.ir\hspace{1cm}sanati@math.um.ac.ir}}
\date{ }
\maketitle

\thispagestyle{myheadings}

\begin{abstract}
 In this paper, we are going to look at the $c$-nilpotent multiplier of a group $G$,
${\cal N}_cM(G)$, as a functor from the category of all groups, ${\cal G}roup$, to the category
of all abelian groups, ${\cal A}b$, and focusing on some functional properties of it. In fact,
by using some results of the first author and others and finding an explicit formula for the
$c$-nilpotent multiplier of a finitely generated abelian group, we try to concentrate on the
commutativity of the above functor with the two famous functors Ext and Tor.\\

{\it A.M.S. Classification 2000}: 20E10, 20K40.
\end{abstract}

\begin{center}
{\bf 1. Introduction}
\end{center}

 Let $G\cong F/R$ be a group, presented as a quotient group of a free group $F$ by a normal
subgroup $R$. Then the {\it Baer-invariant} of $G$, after R. Baer
[1], with respect to the variety $\cal V$, denoted by ${\cal
V}M(G)$, is defined to be
$${\cal V}M(G)=\frac {\displaystyle R\cap V(F)}{\displaystyle[RV^*F]}\ ,$$
where $V(F)$ is the verbal subgroup of $F$ with respect to $\cal
V$ and
$$[RV^*F]=<v(f_1,\ldots ,f_{i-1},f_ir,f_{i+1},\ldots ,f_n)v(f_1,\ldots
,f_i,\ldots f_n)^{-1}\ |\ r\in R,$$
$$1\leq i\leq n,v\in
V,f_i\in F,n\in {\bf N}>.$$

 It can be proved that the Baer-invariant of a group $G$ is independent of the choice of the
presentation of $G$ and it is always an abelian group (See [8]).

 In particular, if $\cal V$ is the variety of abelian groups, $\cal A$, then the Baer-invariant
of $G$ will be $(R\cap F')/[R,F]$, which, following Hopf [6], is
isomorphic to the second cohomology group of $G$, $H_2(G,{\bf
C}^*)$, in finite case, and also is isomorphic to the well-known
notion the {\it Schur multiplier} of $G$, denoted by $M(G)$. The
multiplier $M(G)$ arose in Schur's work [15] of 1904 on
projective representations of a group, and has subsequently found
a variety of other applications. The survey article of Wiegold
[19] and the books by Beyl and Tappe [2] and Karpilovsky [7] form
a fairly comprehensive account of $M(G)$.

 If ${\cal V}$ is the variety of nilpotent groups of class at
most $c\geq 1$, ${\cal N}_c$, then the Baer-invariant of the group $G$ will be
$${\cal N}_cM(G)=\frac {\displaystyle R\cap \ga_{c+1}(F)}{\displaystyle[R,\ _cF]}\ ,$$
where ${\ga}_{c+1}(F)$ is the $(c+1)$st term of the lower central
series of $F$ and $[R, _1F]=[R,F]\ ,\ [R, _cF]=[[R, _{c-1}F],F]$,
inductively. The above notion is also called the $c$-nilpotent
multiplier of $G$ and denoted by $M^{(c)}(G)$ (see [3]).

 The following theorem permit us to look at the notion of the Baer-invariant as a
 functor.\\
{\bf Theorem 1.1}.

 Let $\cal V$ be an arbitrary variety of groups.
Then, using the notion of the Baer-invariant, we can consider the
following covariant functor from the category of all groups,
${\cal G}roup$, to the category of all abelian groups, ${\cal A}b$
$${\cal V}M(-)\ :\ {\cal G}roup\lra {\cal A}b\ ,$$
which assigns to any group $G$ the abelian group ${\cal V}M(G)$.\\
{\bf Proof.} Let $G$ be an arbitrary group. By the properties of
the Baer-invariant, ${\cal V}M(G)$ is independent of the choice
of a presentation of $G$ and it is always abelian. So ${\cal
V}M(-)$ assigns an abelian group to each group $G$. Also, if
$G_1$ and $G_2$ are two arbitrary groups with the following
presentations:
$$ 1\lra R_1\lra F_1\stackrel{\pi_1}{\lra}G_1\lra 1\ \ \
,\ \ \ 1\lra R_2\lra F_2\stackrel{\pi_2}{\lra}G_2\lra 1\ , $$
 and if $\phi :G_1\ra G_2$ is a homomorphism, then, using the universal
property of free groups, there exists a homomorphism $\ol
{\phi}:F_1\ra F_2$. It is easy to see that $\ol {\phi}$ induces a
homomorphism
$$\hat {\phi}:\frac{\displaystyle R_1\cap
V(F_1)}{\displaystyle[R_1V^*F_1]}\lra \frac{\displaystyle R_2\cap
V(F_2)}{\displaystyle[R_2V^*F_2]}\ ,$$
i.e. $\ol{\phi}:{\cal
V}M(G_1)\lra {\cal V}M(G_2)$ is a homomorphism from the
Baer-invariant of $G_1$ to the Baer-invariant of $G_2$. It is a
routine verification to see that the above assignment is a functor from ${\cal G}roup$ to
${\cal A}b$ (see also [8]).$\Box$\\
\newpage
\begin{center}
{\bf $\S$2. Elementary Results}
\end{center}

 Being additive is usually one of the important property that a
functor may have. Unfortunately, the $c$-nilpotent multiplier
functor ${\cal N}_cM(-)$ is {\it not} additive even if we restrict
ourself to abelian groups. The following theorems can prove this
claim.\\
{\bf Theorem 2.1} (I. Schur [14], J. Wiegold [16]).

 Let $G=A\times B$ be the direct product of two groups $A$ and
$B$. Then $$ M(G)\cong M(A)\oplus M(B)\oplus (A_{ab}\otimes
B_{ab})\ .$$ {\bf Theorem 2.2} (B. Mashayekhy and M.R.R.
Moghaddam [11]).

 Let $G\cong {\bf Z}_{n_1}\oplus {\bf Z}_{n_2}\oplus \ldots \oplus {\bf Z}_{n_k}$, be a finite abelian group, where $n_{i+1}|n_i$ for all $1\leq i\leq k-1$ and $k\geq 2$. Then, for all $c\geq 1$, the $c$-nilpotent multiplier of $G$ is
$$ {\cal N}_cM(G)\cong {\bf Z}_{n_2}^{(b_2)}\oplus {\bf
Z}_{n_3}^{(b_3-b_2)}\oplus \ldots \oplus {\bf
Z}_{n_k}^{(b_k-b_{k-1})} \ \ \ \ ,$$ where ${\bf Z}_m^{(n)}$
denotes the direct sum of $n$ copies of the cyclic group ${\bf
Z}_m$, and $b_i$ is the number of basic commutators of weight
$c+1$ on $i$ letters (see [5]).

 One of the interesting corollary of Theorem 2.2 is that the $c$-nilpotent multiplier functors
 can preserve every elementary abelian $p$-group.\\
{\bf Corollary 2.3}.

 Let $G={\bf Z}_p\oplus \ldots \oplus {\bf Z}_p$ ($k$-copies) be
an elementary abelian $p$-group. Then, for all $c\geq 1$, ${\cal
N}_cM(G)$ is also an elementary abelian $p$-group.\\
{\bf Proof.} By Theorem 2.2 we have $${\cal N}_cM(G)\cong {\bf
Z}_p^{(b_2)}\oplus {\bf Z}_p^{(b_3-b_2)}\oplus \ldots \oplus {\bf
Z}_p^{(b_k-b_{k-1})}={\bf Z}_p\oplus \ldots \oplus {\bf Z}_p\ \
(b_k-copies).$$ Hence the result holds. Note that $|G|=p^n$ and
$|{\cal N}_cM(G)|=p^{b_k} .\Box$

 In 1952, C. Miller [12] proved that the Schur multiplier of a free product is isomorphic to the
 direct sum of the Schur multipliers of the free factors. In other words, he proved that the
 Schur multiplier functor $M(-)$ {\it is coproduct-preserving}.\\
{\bf Theorem 2.4} (C. Miller [12]).

 For any group $G_1$ and $G_2$,
$$ M(G_1*G_2)\cong M(G_1)\oplus M(G_2)\ ,$$
where $G_1*G_2$ is the free product of $G_1$ and $G_2$.

 Now, with regards to the above theorem, it seems natural to ask whether the $c$-nilpotent
multiplier functors, ${\cal N}_cM(-),\ c\geq 2$, are
coproduct-preserving or not. To answer the question, first we
state an important theorem of J. Burns and G. Ellis [3,
Proposition 2.13 and its Erratum] which is proved by a
homological method.\\
{\bf Theorem 2.5} (J. Burns and G. Ellis [3]).

 Let $G$ and $H$ be two arbitrary groups, then there is an isomorphism
$${\cal N}_2M(G*H)\cong$$ $$ {\cal N}_2M(G)\oplus {\cal
N}_2M(H)\oplus (M(G)\otimes H_{ab})\oplus (G_{ab}\otimes
M(H))\oplus Tor_1^{{\bf Z}}(G_{ab},H_{ab}) \ .$$

 Now, using the above theorem and properties of tensor product and $Tor_1^{{\bf Z}}$, we can
prove that the second nilpotent multiplier functor ${\cal
N}_2M(-)$, preserves the coproduct of a finite family of cyclic
groups of mutually coprime order.\\
{\bf Corollary 2.6}.

Let $\{{\bf Z}_{n_i}|1\leq i\leq m\}$ be a family of cyclic
groups of mutually coprime order. Then
$${\cal
N}_2M(\prod_{i=1}^{m}\!^{*}{\bf Z}_{n_i})\cong \oplus
\sum_{i=1}^{m}{\cal N}_2M({\bf Z}_{n_i})\ ,$$ where
$\prod_{i=1}^{m}\!^{*}{\bf Z}_{n_i}$ is the free product of ${\bf
Z}_{n_i}$'s, $1\leq i\leq n$.\\
{\bf Proof.} By using  induction on $m$ and the following
properties the result holds.
$${\cal N}_2M({\bf Z}_{n_i})\cong 1, Tor_1^{{\bf Z}}({\bf Z}_{n_i},{\bf Z}_{n_j})\cong{\bf Z}_{n_i}\otimes{\bf Z}_{n_j}=1,\ for\ all\ i\neq j. \Box$$

 Note that the first author has generalized the above corollary to the variety of nilpotent
groups of class at most $c, {\cal N}_c$, for all $c\geq 2$ as follows.\\
{\bf Theorem 2.7} (B. Mashayekhy [10]).

Let $\{{\bf Z}_{n_i}|1\leq i\leq m\}$ be a family of cyclic groups of mutually coprime
order. Then
$${\cal N}_cM(\prod_{i=1}^{m}\!^{*}{\bf Z}_{n_i})\cong
\oplus\sum_{i=1}^{m}{\cal N}_cM({\bf Z}_{n_i})\ ,\ for\ all\ c\geq 1.$$

 In the following example, we are going to show that the condition of being mutually coprime
order for the family of cyclic groups $\{{\bf Z}_{n_i}|1\leq i\leq m\}$ is very essential in the above
results. In other words, we show that the second nilpotent multiplier functor, ${\cal N}_2M(-)$,
 {\it is not coproduct preserving}, in general.\\
{\bf Example}.

 Let $D_{\infty}=<a,b|a^2=b^2=1>\cong {\bf Z}_2*{\bf Z}_2$ be the infinite
dihedral group. Then
$${\cal N}_2M(D_{\infty})\not\cong {\cal N}_2M({\bf Z}_2)\oplus {\cal N}_2M({\bf Z}_2)\
\ .$$
{\bf Proof.} By Theorem 2.5 we have\\
${\cal N}_2M(D_{\infty})\cong {\cal N}_2M({\bf Z}_2)\oplus
{\cal N}_2M({\bf Z}_2) \oplus {\bf Z}_2\otimes M({\bf Z}_2)\oplus
M({\bf Z}_2)\otimes {\bf Z}_2\oplus Tor_1^{{\bf Z}}({\bf Z}_2,{\bf Z}_2)$

$\hspace{1.5cm}\cong Tor_1^{{\bf Z}}({\bf Z}_2,{\bf Z}_2)\cong {\bf Z}_2\otimes {\bf
Z}_2\cong {\bf Z}_2\ .$\\
 But ${\cal N}_2M({\bf Z}_2)\oplus {\cal N}_2M({\bf Z}_2)=1$. Hence the result holds. $\Box$\\
{\bf Note}.

 In 1980 M.R.R. Moghaddam [12] proved that in general, the
Baer-invariant functor commutes with direct limit of a directed
system of groups.\\

 We know that every functor can preserve any
split exact sequence as a split sequence. This property gives us
the following interesting result.\\ \ \ \\
\ \ \\ {\bf Theorem 2.8}.

 Let $G=T\stackrel {\theta}{\rhd\!\!\!<}N$ be the semidirect product
(splitting extension) of $N$ by $T$ under $\theta $. Then ${\cal
V}M(T)$ is a direct summand of ${\cal V}M(G)$, for every
variety of groups $\cal V$.\\

 Note that K.I. Tahara [15] 1972, and W. Haebich [4] 1977, tried to obtain a result
 similar to
the above theorem for the Schur multiplier of a semidirect
product with an emphasis on finding the structure of the
complementary factor $M(T)$ of $M(G)$, as much as possible. Also,
a generalization of Haebich's result [4] presented by the first
author in [9].\\

 Finally, the properties of right and left exactness are some of
the most interesting properties that a functor may have. In the
following, we show that the $c$-nilpotent multiplier functors
{\it are not right or left exact}.\\
{\bf Theorem 2.9}.

 For every $c\geq 1$, the $c$-nilpotent
multiplier functor, ${\cal N}_cM(-)$, is not right exact.\\
{\bf Proof.} Let $G$ be a group such that ${\cal N}_cM(G)\neq 1$
(note that by Theorem 2.2, we can always find such a group $G$).
Let $F$ be a free group and $\pi :F\ra G$ be an epimorphism (we
can always consider a free presentation for a group $G$). Now by
definition of the Baer-invariant we have ${\cal N}_cM(F)=1$
(consider the free presentation $1\ra 1\ra F\ra F\ra 1$ for $F$).
Therefore, it is easy
to see that ${\cal N}_cM(F)\lra{\cal N}_cM(G)$ is not onto. $\Box$\\
{\bf Theorem 2.10}.

The $c$-nilpotent multiplier functor, ${\cal N}_cM(-)$, is not
left exact, in general.\\
{\bf Proof.} Suppose $G={\bf Z}_4\oplus{\bf Z}_4$. Then by
Theorem 2.1 we have
$$M(G)\cong M({\bf Z}_4)\oplus M({\bf Z}_4)\oplus ({\bf Z}_4\otimes {\bf Z}_4)\cong
{\bf Z}_4\ .$$
By a famous result on the Schur multiplier we know that every finite $p$-group can
be embedded in a finite $p$-group whose Schur multiplier is
elementary abelian $p$-group (see [7,17]). So there exists an
exact sequence $G\stackrel{\theta}{\ra} H\ra 1$, where $H$ is a finite
$2$-group and $M(H)$ is an elementary abelian $2$-group. Hence
$M(\theta ):M(G)\ra M(H)$ can not be a monomorphism. $\Box$

\begin{center}
{\bf 3. Main Results}
\end{center}

 In this section, we will see the behaviour of the functor ${\cal
N}_cM(-)$ with the functors $Ext^n_{\bf Z}({\bf Z_m},-)$ and
$Tor_n^{\bf Z}({\bf Z_m},-)$. First, by using notations and
similar method of paper [11], we can present an explicit formula
for the $c$-nilpotent multiplier of a finitely generated
abelian groups as follows.\\
{\bf Theorem 3.1}.

 Let $G\cong{\bf Z}^{(n)}\oplus{\bf Z}_{n_1}\oplus {\bf Z}_{n_2}\oplus \ldots \oplus
{\bf Z}_{n_k}$, be a finitely generated abelian group, where
$n\geq 0$, $n_{i+1}|n_i$ for all $1\leq i\leq k-1$ and $k\geq 2$.
Then, for all $c\geq 1$, the $c$-nilpotent multiplier of $G$ is
$${\cal N}_cM(G)\cong{\bf Z}^{(b_n)}\oplus{\bf
Z}_{n_1}^{(b_{n+1}-b_n)}\oplus\ldots\oplus{\bf
Z}_{n_k}^{(b_{n+k}-b_{n+k-1})} \ \ \ \ ,$$
where $b_1=b_0=0$\\
{\bf Proof.} Clearly ${\bf Z}\otimes{\bf Z}\cong{\bf Z}$, ${\bf
Z}\otimes{\bf Z}_{n_i} \cong{\bf Z}_{n_i}$ and ${\bf
Z}_{n_i}\otimes{\bf Z}_{n_{i+1}}\cong{\bf Z}_{n_{i+1}}$. Hence we
have
\begin{center}
${\bf Z}^{(t)}\otimes{\bf Z}_{n_1}\otimes{\bf
Z}_{n_2}\otimes\ldots\otimes{\bf
Z}_{n_r}\stackrel{(*)}{\cong}{\bf Z}_{n_r}$ \ \ and\ \ ${\bf
Z}\otimes\ldots\otimes{\bf Z}\cong{\bf Z}.$
\end{center}
for all $t\geq 0$
and $r\geq 1$. Thus by theorem 2.3 of [11] we have
$${\cal N}_cM({\bf Z}^{(n)})\cong T({\bf Z},\ldots,{\bf
Z})_{c+1}\cong{\bf Z}^{(b_n)}.$$ We remind that
$T(H_1,\ldots,H_n)_{c+1}$ is the summation of all the tensor
products corresponding to the subgroup generated by all the basic
commutators of weight $c+1$ on $n$ letters $x_1,\ldots,x_n,$
where $x_i\in H_i$ for all $1\leq i\leq n$. Now, by induction
hypothesis assume
$${\cal N}_cM({\bf Z}^{(n)}\oplus{\bf Z}_{n_1}\oplus{\bf Z}_{n_2}\oplus\ldots\oplus
{\bf Z}_{n_{k-1}})\cong{\bf Z}^{(b_n)}\oplus{\bf
Z}_{n_1}^{(b_{n+1}-b_n)}\oplus\ldots \oplus{\bf
Z}_{n_{k-1}}^{(b_{n+k-1}-b_{n+k-2})} \ .$$ Then we have
$${\cal N}_cM({\bf Z}^{(n)}\oplus{\bf Z}_{n_1}\oplus{\bf
Z}_{n_2}\oplus\ldots\oplus{\bf Z}_{n_k})\cong T(\underbrace{{\bf
Z},\ldots,{\bf Z}}_{n-times},{\bf Z}_{n_1},\ldots,{\bf
Z}_{n_k})_{c+1}$$
$$\hspace{6.cm}\cong T(\underbrace{{\bf Z},\ldots,{\bf Z}}_{n-times},{\bf Z}_{n_1},
\ldots,{\bf Z}_{n_{k-1}})_{c+1} \oplus L\ ,$$ where $L$ is the
summation of all the tensor products of ${\bf Z},{\bf Z}_{n_1},
\ldots,{\bf Z}_{n_k}$ corresponding to the subgroup generated by
all the basic commutators of weight $c+1$ on $n+k$ letters which
involve ${\bf Z}_{n_k}$. Using $(*)$, all those tensor product
are isomorphic to ${\bf Z}_{n_k}$. So $L$ is the direct summand
of $(b_{n+k}-b_{n+k-1})$-copies of ${\bf Z}_{n_k}$. Hence the
result follows by induction. $\Box$

 For the rest of the paper we need the following lemmas.\\
{\bf Lemma 3.2}.

 For any abelian groups $A$ and $B$, we have\\
$(i)$ $Ext^1_{\bf Z}({\bf Z}/m{\bf Z},B)\cong B/mB .$\\
Also, $Ext^n_{\bf Z}(A,B)=0$, for all $n\geq 2$.\\
$(ii)$ If $A$ and $B$ are finite abelian groups, then
$$Ext^1_{\bf Z}(A,B)\cong Ext^1_{\bf Z}(B,A).$$
$(iii)$ $Tor_1^{\bf Z}({\bf Z}/m{\bf Z},B)\cong B[m],$ where
$B[m]=\{b\in B: mb=0\}$. Also, $Tor_n^{\bf Z}(A,B)=0,$ for all
$n\geq 2$, and  $Tor_1^{\bf Z}(A,B)\cong Tor_1^{\bf Z}(B,A)$.\\
{\bf Proof.} See [14, Chapters 7, 8]. $\Box$\\ \ \ \\
{\bf Lemma 3.3}.

Let $A$ and $\{B_k\}_{k\in I}$ be abelian groups. Then for all $n\geq 0$ the
following isomorphism hold.\\
$$(i)\ Ext^n_{\bf Z}(A,\prod_{k\in I}B_k)\cong \prod_{k\in I}
Ext^n_{\bf Z}(A,B_k),\ Ext^n_{\bf Z}(\coprod_{k\in I}B_k,A)\cong
\prod_{k\in I}Ext^n_{ \bf Z}(B_k,A).$$
 $$(ii)\ Tor_n^{\bf
Z}(A,\coprod_{k\in I}B_k)\cong\coprod_{k\in I}Tor_n^{\bf
Z}(A,B_k),\ Tor_n^{\bf Z}(\coprod_{k\in
I}B_k,A)\cong\coprod_{k\in I}Tor_n^{\bf Z}(B_k,A).$$
 {\bf Proof.} See [14]. $\Box$

It is obvious that the functor ${\cal N}_cM(-)$ commutes with the
functors $Ext^n_{\bf Z}({\bf Z}_m,-)$, and $Tor_n^{\bf Z}({\bf
Z}_m,-)$ for all $n\geq 2$, by lemma 3.2. Now we are going to pay
our attention to the functors
$Ext^1_{\bf Z}({\bf Z}_m,-)$, $Ext^1_{\bf Z}(-,{\bf Z}_m)$, and $Tor_1^{\bf Z}({\bf Z}_m,-)$.\\
{\bf Theorem 3.4}.

Let $D\cong{\bf Z}^{(n)}\oplus{\bf Z}_{n_1}\oplus {\bf Z}_{n_2}\oplus \ldots \oplus
{\bf Z}_{n_k}$, be a finitely generated abelian group, where
$n\geq 0$, $n_{i+1}|n_i$ for all $1\leq i\leq k-1$. Then, for all
$c\geq 1$, the following isomorphisms hold.\\
$(i)$  ${\cal N}_cM{\big (}Ext^1_{\bf Z}({\bf Z}_m,D){\big
)}\cong{\bf Z}_m^{(b_n)} \oplus (\oplus\sum_{i=1}^k{\bf
Z}_{(n_i,m)}^{(b_{n+i}-b_{n+i-1})}).$\\
$(ii)$  $Ext^1_{\bf Z}({\bf Z}_m,{\cal N}_cM(D))\cong{\bf
Z}_m^{(b_n)}\oplus (\oplus\sum_{i=1}^k{\bf
Z}_{(n_i,m)}^{(b_{n+i}-b_{n+i-1})}).$\\
$(iii)$  ${\cal N}_cM(Ext^1_{\bf Z}(D,{\bf Z}_m))\cong \oplus
\sum_{i=2}^k{\bf Z}_{(n_i,m)}^{(b_{i}-b_{i-1})}.$\\
$(iv)$  $Ext^1_{\bf Z}({\cal N}_cM(D),{\bf Z}_m)\cong
\oplus\sum_{i=1}^k{\bf
Z}_{(n_i,m)}^{(b_{n+i}-b_{n+i-1})}.$\\
$(v)$  ${\cal N}_cM(Tor_1^{\bf Z}({\bf
Z}_m,D))\cong\oplus\sum_{i=2}^k{\bf
Z}_{(n_i,m)}^{(b_i-b_{i-1})}.$\\
$(vi)$  $Tor_1^{\bf Z}({\bf Z}_m,{\cal
N}_cM(D))\cong\oplus\sum_{i=1}^k{\bf
Z}_{(n_i,m)}^{(b_{n+i}-b_{n+i-1})}.$\\
{\bf Proof.} $(i)$ By Lemma 3.3(i), $Ext^1_{\bf Z}({\bf Z}/m{\bf
Z},{\bf Z})\cong {\bf Z}/m{\bf Z}\cong{\bf Z}_m$. Now by using
Lemmas 3.3(i) and 3.2(i), we have  $$Ext^1_{\bf Z}({\bf
Z}_m,D)\cong(Ext^1_{\bf Z}({\bf Z}_m,{\bf
Z}))^{(n)}\oplus(\oplus\sum_{i=1}^kExt^1_{\bf Z}({\bf Z}_m,{\bf
Z}_{n_i}))$$
$$\cong{\bf Z}_m^{(n)}\oplus(\oplus\sum_{i=1}^k{\bf Z}_{n_i}/m{\bf Z}_{n_i}).$$
One can see that for every $n,m\in {\bf Z}$, we have ${\bf
Z}_m/n{\bf Z}_m\cong {\bf Z}_{(n,m)}$. Therefore
$$Ext^1_{\bf Z}({\bf Z}_m,D)\cong{\bf Z}_m^{(n)}\oplus(\oplus\sum_{i=1}^k{\bf Z}_{
(n_i,m)}).$$ Now, by Theorem 2.2 and by noting that
$(m,n_{i+1})|(m,n_i)|m$ we have
$${\cal N}_cM(Ext^1_{\bf Z}({\bf Z}_m,D))$$
$$\cong{\bf Z}_m^{(b_2-b_1)}\oplus{\bf Z}_m^{(b_3-b_2)}\oplus\ldots\oplus{\bf Z}_m^{
(b_n-b_{n-1})}\oplus{\bf
Z}_{(n_1,m)}^{(b_{n+1}-b_{n})}\oplus\ldots\oplus{\bf Z}_{
(n_k,m)}^{(b_{n+k}-b_{n+k-1})}$$
 $$\cong{\bf Z}_m^{(b_n)}\oplus(\oplus\sum_{i=1}^k{\bf
Z}_{(n_i,m)}^{(b_{n+i}-b_{n+i-1})}).$$
$(ii)$ By Theorem 3.1 and
Lemmas 3.3(i) and 3.2(i), we have
\begin{center}
$Ext^1_{\bf Z}({\bf Z}_m,{\cal N}_cM(D))\cong Ext^1_{\bf Z}({\bf Z}_m,{\bf Z})^{
(b_n)}\oplus(\oplus\sum_{i=1}^k(Ext^1_{\bf Z}({\bf Z}_m,{\bf
Z}_{n_i}))^{ (b_{n+i}-b_{n+i-1})})$

$\cong{\bf Z}_m^{(b_n)}\oplus (\oplus\sum_{i=1}^k{\bf
Z}_{(n_i,m)}^{(b_{n+i}-b_{n+i-1})}).$
\end{center}
$(iii)$ By Lemmas 3.3(ii) and 3.2(ii) we have
$$Tor_1^{\bf Z}({\bf Z}_m,D)\cong(Tor_1^{\bf Z}({\bf Z}_m,{\bf
Z}))^{(n)}\oplus(\oplus\sum_{i=1}^kTor_1^{\bf Z}({\bf Z}_m,{\bf
Z}_{n_i})) \cong\oplus\sum_{i=1}^k{\bf Z}_{n_i}[m].$$ Note that
$Tor_1^{\bf Z}({\bf Z}_m,{\bf Z})\cong 1$ and ${\bf
Z}_{n}[m]\cong{\bf Z}_{(m,n)}$. So we have $Tor_1^{\bf Z}({\bf
Z}_m,D)\cong\oplus\sum_{i=1}^k{\bf Z}_{(n_i,m)}$. Now by Theorem
2.2 the result holds.\\
$(iv)$ Again by using Theorem 3.1 and Lemmas 3.3(ii) and 3.2(ii),
we have
$$Tor_1^{\bf Z}({\bf Z}_m,{\cal N}_cM(D))\cong(Tor_1^{\bf
Z}({\bf Z}_m, {\bf Z}))^{(b_n)}\oplus(\oplus\sum_{i=1}^kTor_1^{\bf
Z}({\bf Z}_m,{\bf Z}_{n_i}^{(b_{n+i}-b_{n+i-1})})$$
$$\cong\oplus\sum_{i=1}^kTor_1^{\bf Z}({\bf Z}_m,{\bf
Z}_{n_i}^{(b_{n+i}-b_{n+i-1})})\cong \oplus\sum_{i=1}^k{\bf
Z}_{(n_i,m)}^{(b_{n+i}-b_{n+i-1})}. \ \Box$$

 In the following corollary you can find some of main results of
 the paper.\\
{\bf Corollary 3.5}.

 Let $D$ be an arbitrary finitely generated abelian group. Then\\
$(i)$  ${\cal N}_cM{\big (}Ext^1_{\bf Z}({\bf Z}_m,D){\big )}\cong
Ext^1_{\bf Z}({\bf Z}_m,{\cal N}_cM(D)).$\\
$(ii)$ If $D$ is also finite, then
$${\cal N}_cM(Tor_1^{\bf Z}({\bf Z}_m,D))\cong Tor_1^{\bf Z}({\bf Z}_m,{\cal N}_cM(D)),$$\\
$${\cal N}_cM(Ext^1_{\bf Z}(D,{\bf Z}_m))\cong Ext^1_{\bf Z}({\cal N}_cM(D),{\bf Z}_m).$$
$(iii)$ If $D$ is infinite, then
$${\cal N}_cM(Tor_1^{\bf Z}({\bf Z}_m,D))\not\cong Tor_1^{\bf Z}({\bf Z}_m,
{\cal N}_cM(D)).$$

$${\cal N}_cM(Ext^1_{\bf Z}(D,{\bf Z}_m))\not\cong Ext^1_{\bf Z}({\cal N}_cM(D),
{\bf Z}_m).$$ This means that the $c$-nilpotent multiplier
functors, ${\cal N}_cM(-)$ {\it do not} commute with $Tor_1^{\bf
Z}({\bf Z}_m, -)$ and $Ext^1_{\bf Z}(-,{\bf Z}_m)$, in infinite case.\\
{\bf Proof.} $(i)$ It is clear by parts $(i),\ (ii)$ of the
previous theorem.\\
$(ii)$ By putting $n=0$ in parts $(iii)$ to $(vi)$ of the
previous theorem, the result holds.\\
$(iii)$ Since $D$ is infinite, so $n\geq 1$. Hence the result
holds by the previous theorem parts $(iii)$ to $(vi)$.\ $\Box$\\

 We know that $Hom({\bf Z}_m,{\bf Z})\cong 0$ and $Hom({\bf Z},{\bf Z}_m)\cong{\bf Z}_m$.
So by similar methods of Theorem 3.4 we are going to indicate the
behaviour of functor ${\cal N}_cM(-)$ with $Ext_{\bf Z}^0({\bf
Z}_m,-)=Hom({\bf Z}-m,-)$, $Ext_{\bf Z}^0(-,{\bf Z}_m)=Hom(-,{\bf
Z}_m)$, and
$Tor^{\bf Z}_0({\bf Z}_m,-)={\bf Z}_m\otimes -$ as the following theorem.\\
{\bf Theorem 3.6}.

 For any finitely generated abelian group $D\cong{\bf Z}^{(n)}\oplus{\bf Z}_{n_1}
\oplus {\bf Z}_{n_2}\oplus \ldots \oplus {\bf Z}_{n_k}$, we have\\
$(i)$ ${\cal N}_cM(Hom({\bf Z}_m,D)\cong{\bf
Z}_{(m,n_2)}^{(b_2)}\oplus\ldots \oplus{\bf
Z}_{(m,n_k)}^{(b_{k}-b_{k-1})}.$\\
$(ii)$ $Hom({\bf Z}_m,{\cal N}_cM(D))\cong{\bf
Z}_{(m,n_1)}^{(b_{n+1}-b_n)}\oplus \ldots\oplus{\bf
Z}_{(m,n_k)}^{(b_{n+k}-b_{n+k-1})}.$\\
$(iii)$ If $D$ is finite, then ${\cal N}_cM(Hom({\bf Z}_m,D)\cong
Hom({\bf Z}_m,{\cal N}_cM(D)).$\\
If $D$ is infinite, then ${\cal N}_cM(Hom({\bf Z}_m,D)\not\cong
Hom({\bf Z}_m,{\cal N}_cM(D)).$\\
$(iv)\ {\cal N}_cM(Hom(D,{\bf Z}_m))\cong Hom({\cal N}_cM(D),{\bf
Z}_m)$$ $$\cong {\bf Z}_m^{(b_n)}\oplus{\bf
Z}_{(m,n_1)}^{(b_{n+1}-b_n)}\oplus\ldots\oplus {\bf
Z}_{(m,n_k)}^{(b_{n+k}-b_{n+k-1})}.$\\
$(v)\ {\cal N}_cM({\bf Z}_m\otimes D)\cong{\bf Z}_m\otimes{\cal
N}_cM(D)$$ $$\cong{\bf Z}_m^{(b_n)}\oplus{\bf
Z}_{(m,n_1)}^{(b_{n+1}-b_n)}\oplus\ldots\oplus {\bf
Z}_{(m,n_k)}^{(b_{n+k}-b_{n+k-1})}.$

 Now, in the following we  are going to show that our conditions in
the previous results are essential. In general case $Ext^i_{\bf
Z}(A,-)$ and $Tor_i^{\bf Z}(A,-)$,
where $A$ is not cyclic, {\it do not commute} with ${\cal N}_cM(-)$, for $i=0,1$.\\
{\bf Some Examples}.

$(a)$ ${\cal N}_cM(Ext^1_{\bf Z}({\bf Z}_n\oplus{\bf Z}_n,{\bf
Z}_n))\cong{\bf Z}_n^{(b_2)}\not\cong 1\cong Ext^1_{\bf Z}({\bf
Z}_n\oplus{\bf Z}_n,{\cal N}_cM({\bf Z}_n)),$ i.e
$${\cal N}_cM(Ext_{\bf Z}^1(-,A))\not\cong Ext_{\bf Z}^1({\cal N}_cM(-),A).$$
$(b)$ ${\cal N}_cM(Ext^1_{\bf Z}({\bf Z}_n,{\bf Z}_n\oplus{\bf
Z}_n)\cong{\bf Z}_n^{(b_2)}\not\cong 1\cong Ext^1_{\bf Z}({\cal
N}_cM({\bf Z}_n),{\bf Z}_n\oplus{\bf Z}_n),$ i.e
$${\cal N}_cM(Ext_{\bf Z}^1(A,-))\not\cong Ext_{\bf Z}^1(A,{\cal N}_cM(-)).$$
$(c)$ ${\cal N}_cM(Tor_1^{\bf Z}({\bf Z}_n\oplus{\bf Z}_n,{\bf
Z}_n))\cong{\bf Z}_n^{(b_2)}\not\cong 1\cong Tor_1^{\bf Z}({\bf
Z}_n\oplus{\bf Z}_n,{\cal N}_cM({\bf Z}_n)),$ i.e
$${\cal N}_cM(Tor^{\bf Z}_1(-,A))\not\cong Tor^{\bf Z}_1({\cal N}_cM(-),A).$$
$(d)$ ${\cal N}_cM(({\bf Z}_n\oplus{\bf Z}_n\otimes{\bf
Z}_n))\cong{\bf Z}_n^{(b_2)}\not\cong ({\bf Z}_n\oplus{\bf
Z}_n\otimes{\cal N}_cM({\bf Z}_n)),$ i.e
$${\cal N}_cM(A\otimes -))\not\cong (A\otimes{\cal N}_cM(-)).$$
$(e)$ ${\cal N}_cM(Hom({\bf Z}_n\oplus{\bf Z}_n,{\bf Z}_n)\cong
{\bf Z}_n^{(b_2)}\not\cong 1\cong Hom({\bf Z}_n\oplus{\bf
Z}_n,{\cal N}_cM({\bf Z}_n)),$ i.e.
$${\cal N}_cM(Hom(A,-))\not\cong Hom(A,{\cal N}_cM(-)).$$
$(f)$ ${\cal N}_cM(Hom({\bf Z}_{14}\oplus{\bf Z}_2,{\bf
Z}_6\oplus{\bf Z}_3))\cong {\bf Z}_2^{(b_2)} \not\cong 1\cong
Hom({\bf Z}_14\oplus {\bf Z}_2,{\bf Z}_3^{(b_2)}\cong Hom({\bf
Z}_14\oplus {\bf Z}_2,{\cal N}_cM({\bf Z}_6\oplus{\bf Z}_3)),$
i.e.
$${\cal N}_cM(Hom(A,-))\not\cong Hom(A,{\cal N}_cM(-)).$$
$(g)$ ${\cal N}_cM(Hom({\bf Z}_6\oplus{\bf Z}_2,{\bf
Z}_9\oplus{\bf Z}_3))\cong{\bf Z}_3^{(b_2)}\not\cong 1\cong
Hom({\bf Z}_2^{(b_2)},{\bf Z}_9\oplus{\bf Z}_3))\cong Hom({\cal
N}_cM({\bf Z}_6\oplus{\bf Z}_2),{\bf Z}_9\oplus{\bf Z}_3),$ i.e
$${\cal N}_cM(Hom(-,A))\not\cong Hom({\cal N}_cM(-),A).$$
$(h)$ $M(Hom(D,{\bf Z}_m))\not\cong Hom(M(D),{\bf Z}_m)$, and
$M(D\otimes{\bf Z}_m)\not\cong M(D)\otimes{\bf Z}_m,$\\
when $D$ is not abelian; Because one can see that $Hom(S_n,{\bf
Z}_2)\cong{\bf Z}_2$, for $n\geq 2$. Also we know that
$M(S_n)\cong{\bf Z}_2$, for each $n\geq 4$, see [7, theorem
2.12.3]. Now

$$1\cong M(Hom(S_n,{\bf Z}_2))\not\cong Hom(M(S_n),{\bf Z}_2)\cong{\bf
Z}_2,$$ Moreover $S_n\otimes{\bf Z}_2\cong S_n/S_n'\otimes{\bf
Z}_2\cong{\bf Z}_2\otimes{\bf Z}_2\cong {\bf Z}_2.$ Then
$$1\cong M(S_n\otimes{\bf Z}_2))\not\cong M(S_n)\otimes{\bf Z}_2\cong{\bf
Z}_2.$$ The functor ${\cal S}=A\otimes -$, where $A$ is a
non-cyclic group {\it does not commute} with the fuctor $ {\cal
N}_cM(-)$. Put $A={\bf Z}_{n_1}\oplus{\bf Z}_{n_2}, G={\bf Z}_n$,
where $n_2|n_1$. Then $A\otimes G\cong{\bf Z}_{m_1}\oplus {\bf
Z}_{m_2}$, where $m_i=(n,n_i)$, for $i=1,2$. Clearly $m_2|m_1$,
so by Theorem 3.1 we have ${\cal N}_cM(A\otimes G)\cong{\bf
Z}_{m_2}^{(b_2)}$. On the other hand, we have $A\otimes{\cal
N}_cM(G)\cong A\otimes 1=1$. Hence ${\cal N}_cM(A\otimes
G)\not\cong A\otimes{\cal N}_cM(G)$.

We should also point out that the Theorem 3.1 shows that the
$c$-nilpotent multiplier functor, ${\cal N}_cM(-)$, {\it does not
preserve} the tensor product, for
\begin{center}
${\cal N}_cM({\bf Z}_m\otimes G_{ab})\not\cong{\cal N}_cM({\bf Z}_m)\otimes{\cal N}_c
M(G_{ab})=1\ .\ \Box$
\end{center}

\end{document}